\documentclass{article}   	
\usepackage{geometry}       
\usepackage{graphicx}
\usepackage[usenames,dvipsnames]{xcolor}
\usepackage{amssymb,amsmath,amsthm,amsfonts}

\allowdisplaybreaks

\usepackage[english]{babel}
\usepackage[utf8]{inputenc}
\usepackage[colorlinks]{hyperref}
\hypersetup{linkcolor=blue,citecolor=blue,filecolor=black,urlcolor=blue}

\usepackage{mathtools}
\usepackage{dsfont}

\usepackage[sc]{mathpazo}

\usepackage{tikz}
\newtheorem{proposition}{Proposition}[section]
\newtheorem{theorem}[proposition]{Theorem}
\newtheorem{corollary}[proposition]{Corollary}
\newtheorem{lemma}[proposition]{Lemma}

\theoremstyle{definition}
\newtheorem{definition}[proposition]{Definition}
\theoremstyle{remark}
\newtheorem{remark}[proposition]{Remark}
\numberwithin{equation}{section}

\newcommand{\eps}{\varepsilon}

\newcommand{\N}{{\mathbb{N}}}

\newcommand{\R}{{\mathbb{R}}}

\newcommand{\sphere}{{\mathbb{S}}}

\newcommand{\loc}{{\text{loc}}}
\newcommand{\per}{{\text{per}}}

\newcommand{\Ccal}{{\mathcal{C}}}

\newcommand{\Lcal}{{\mathcal{L}}}
\newcommand{\Mcal}{{\mathcal{M}}}

\DeclareMathOperator{\bari}{bar}

\newcounter{fakecounter}

\newcommand{\ind}[1]{\mathds{1}_{#1}}
\newcommand{\labhr}{\mu_1}
\newcommand{\cell}{\Ccal}

\DeclareMathOperator{\od}{\Lambda}

\newcommand{\IM}{\Lambda_0}

\title{Shape optimization of a small favorable region in a periodically fragmented environment}
\author{Gianmaria Verzini}

\begin{document}
\maketitle

\begin{abstract}
We consider a shape  optimization problem for the persistence threshold of a biological 
species dispersing in a periodically fragmented environment, the unknown shape corresponding 
to the portion of the habitat which is favorable to the population. Analytically, this translates 
in the minimization of a weighted eigenvalue of the periodic Laplacian, with respect to a 
bang-bang indefinite weight. For such problem, we exploit some recent results obtained in the 
framework of Dirichlet or Neumann boundary conditions, to provide 
a full description of the singularly perturbed regime in which the volume of the favorable zone 
vanishes.

First, we show that the optimal favorable zone shrinks to a connected, convex, nearly spherical set, 
in $C^{1,1}$ sense. Secondly, we show that the spherical asymmetry of the optimal favorable 
zone decays exponentially, with respect to a negative power of its volume, in the  $C^{1,\alpha}$ sense, 
for every $\alpha<1$. 
\end{abstract}%

\noindent
{\footnotesize \textbf{AMS-Subject Classification}}. 
{\footnotesize 49R05, 49Q10, 92D25, 35P15
}\bigskip\\

\emph{A Sandro, con immensa stima e sincera gratitudine.}

\section{Introduction} \label{sec:intro}

In this note we deal with the asymptotic properties ---both qualitative and quantitative--- of the free boundary 
arising in a shape optimization problem, motivated by studies in population dynamics. We first describe 
in details the motivations from the point of view both of the model and of the previous literature, 
then we turn to the description of our results.

\subsection{Optimal design of a survival threshold} \label{sec:motiv}

The motivation of our study comes from population dynamics. Roughly, we deal with the question: 
"In some environment, determine the optimal spatial arrangement of resources (in particular, 
favorable and unfavorable regions) for a species to survive". In this section we are going to 
describe how this problem translates mathematically into a spectral optimization problem, which in 
turn can be framed as a shape optimization one.

Let us consider a biological species, of density 
$u=u(x,t)$, which disperses in a heterogeneous habitat \cite{MR0043440}, according to the generalized logistic equation
\begin{equation}\label{eq:evlog}
\begin{cases}
u_t - d \Delta u = m(x) u - h(x,u)  & x\in\R^N, t>0,\smallskip\\
u(x,0)=g(x) \gneqq 0 & x\in\R^N.
\end{cases}
\end{equation}

Here $d>0$ is the constant motility coefficient and $g$ is a bounded, non-negative, nontrivial initial 
condition. The reaction term at the right hand side is split into two contributions: in the linear part, 
$m\in L^\infty(\R^N)$ denotes the intrinsic growth rate, which encodes the properties of the habitat, and in particular
$m>0$ a.e. on favorable regions of the habitat, while $m<0$ a.e. on unfavorable ones; on the other hand, 
the positive, superlinear term $h$ accounts for the intraspecies competition, modeled on the prototypical case $h(x,s)=s^2$. 

In this paper we deal with the case of a  \emph{periodically fragmented environment}, meaning 
that both $m$ and $h$ are periodic in each spatial variable $x_i$ (although $g$ needs not). More precisely,
let $(e_1,\dots, e_n)$ denote the canonical orthonormal basis of the euclidean space $\R^N$ and 
$L_1,\dots,L_N$ be fixed positive numbers. We say that a function $f : \R^N \to \R$ is \emph{periodic} 
if $f(x+L_ie_i) = f(x)$ for every $i= 1,\dots,N$ and we denote with $\cell$ the period cell defined by
\[
\cell = \left(-\frac{L_1}{2},\frac{L_1}{2}\right) \times \dots \times \left(-\frac{L_N}{2},
\frac{L_N}{2}\right).
\]
Accordingly, we write 
\[
L^p_{\per}(\cell), \qquad H^1_{\per}(\cell),
\] 
and so on for spaces of periodic functions, with period cell $\cell$.

Problem \eqref{eq:evlog} in a periodic environment, with more general diffusion operators, 
has been considered by Berestycki, Hamel and Roques in a series of papers \cite{MR2214420,MR2155900,haro}. Under mild 
assumptions on $h$ and $g$ one can prove that \eqref{eq:evlog} admits a global-in-time 
solution $u_g$, and one main question concerns the long time behavior of such solution. Roughly, two alternatives 
may occur: either $u_g(\cdot,t)\to0$ as $t\to+\infty$, that is, \emph{extinction} happens; or this is not the case, and 
we talk about \emph{persistence}. Although, in principle, such behavior may depend also on $h$ and $g$, it turns out 
that this is not the case, and the only relevant parameters are $m$ and $d$, via the principal eigenvalue 
$\labhr=\labhr[m,d]$ of the linearized elliptic operator
\[
\Lcal_0\phi := -d \Delta \phi - m \phi
\]
with periodicity conditions. Namely, $\labhr$ is the unique real number such
that there exists a function $\phi_1>0$ which satisfies
\begin{equation}\label{eq:lam_BHR}
\begin{cases}
\Lcal_0 \phi_1 = \labhr \phi \text{ in }\R^N\\
\phi_1\text{ is periodic, }\phi_1>0, 
\end{cases}
\end{equation}
or, alternatively, 
\[
\labhr[m,d]:= \inf \left\{
\int_\cell |\nabla u|^2 - m u^2 \,dx :  
u 
\text{ is periodic, }\int_\cell u^2  \,dx= 1\right\}.
\]
It is a consequence of the results in \cite{MR2214420}, in particular of
Theorems 2.1, 2.4, and 2.6 therein, that:
\begin{enumerate}
\item if $\labhr[m,d] \ge 0$ then $u_g(\cdot,t)\to0$ as $t\to+\infty$, for every $g$ (extinction);
\item if $\labhr[m,d] < 0$ then $u_g(\cdot,t)\to p(\cdot) $ as $t\to+\infty$, for every $g$, 
where $p$ is the unique positive steady state of  \eqref{eq:evlog}, which is periodic (persistence).
\end{enumerate}

Since $\labhr[m,d]$ is increasing in $d$, one can translate these results in a language that 
makes the dependence on $m$ and $d$ more apparent.
\begin{proposition}\label{prop:dstar}
Let $m\in L^\infty(\R^N)$ be periodic. Then
\begin{enumerate}
\item (pointwise hostile habitat) if $m\le0$ a.e. in $\R^N$ then $\labhr[m,d]\ge0$, for every $d>0$;
\item (favorable on average habitat) if 
\[
\int_\cell m  \,dx \ge0\qquad \text{ and }\qquad \operatorname*{ess\,sup}_\cell m > 0
\]
then $\labhr[m,d]<0$, for every $d>0$;
\item\label{it:nontr_case} in the remaining case 
\begin{equation}\label{eq:rem_case}
\int_\cell m  \,dx <0\qquad \text{ and }\qquad \operatorname*{ess\,sup}_\cell m>0
\end{equation}
there exists $d^*=d^*(m)>0$ such that:
\begin{enumerate}
\item if $0<d<d^*(m)$ then $\labhr[m,d]<0$,
\item if $d\ge d^*(m)$ then $\labhr[m,d]\ge0$.
\end{enumerate}
\end{enumerate}
\end{proposition}
This proposition, although not stated therein in these terms, is a straightforward consequence of 
the analysis performed in \cite[Section 5]{MR2214420}. 
We also refer to the book \cite{MR3380662bis}, Chapter 10, which contains 
a detailed analysis of the previously described issues, in the case of Neumann boundary conditions.

In view of Proposition \ref{prop:dstar}, from now on we focus on classes of habitats $m$ which fulfill 
the conditions in case \ref{it:nontr_case}. We obtain that, for every such $m$, 
$d^*(m)$ is characterized by the equation
\[
\labhr[m,d^*(m)] = 0.
\]
In terms of the associated periodic eigenfunction $\phi_1>0$, this can be written as
\begin{equation*}
-d^*(m) \Delta \phi_1 - m \phi_1 = \labhr[m,d^*(m)] \phi_1 = 0,
\qquad
\text{ i.e. }
- \Delta \phi_1 = \frac{1}{d^*(m)} m \phi_1.
\end{equation*}

We deduce that 
\begin{equation}\label{eq:dstar_vs_la}
d^*(m) = \frac{1}{\lambda_1(m)},
\end{equation}
where $\lambda_1=\lambda_1(m)$ is the \emph{positive} principal eigenvalue of the 
weighted problem 
\begin{equation}\label{eq:weighted_eigenproblem}
\begin{cases}
-\Delta u = \lambda m u \text{ in }\R^N\\
u\text{ periodic}, 
\end{cases}
\end{equation}
that is,
\begin{equation}\label{eq:weighted_eigenvalue}
\lambda_1(m):= \inf \left\{
\int_\cell |\nabla u|^2 \,dx :  
u \in H^1_\per(\cell),\ \int_\cell m u^2  \,dx= 1\right\},
\end{equation}
achieved by $\phi_1$. (If $m$ changes sign, problem \eqref{eq:weighted_eigenproblem} has 
two principal eigenvalues; since we are assuming that $m$ has negative average, the other one 
is $\lambda_{-1}(m)=0$, achieved by a constant eigenfunction when the constraint in \eqref{eq:weighted_eigenvalue}
is changed into $\int_\cell m u^2  \,dx= -1$.)

Summarizing, we obtain that, under condition \eqref{eq:rem_case},
\begin{equation}\label{eq:iff_pers}
\text{persistence}
\qquad\iff\qquad
d<\frac{1}{\lambda_1(m)}.
\end{equation}
We deduce that, in this setting, in order to foster the population survival it is convenient 
to decrease $\lambda_1(m)$ as much as possible. 
Accordingly, the question "In some environment, determine the optimal spatial arrangement
of resources (in particular favorable and unfavorable regions) for a species to survive" 
translates mathematically into: "Minimize $\lambda_1(m)$, for $m$ belonging to a suitable class
of weights $\Mcal$". In particular, the most common class of weights considered in the literature 
deals with both pointwise constraints and restrictions of the total amount of resources 
available, such as in
\begin{equation}\label{eq:Mcal}
\Mcal:=\left\{ m \in L^\infty_\per (\cell): 
-\beta \le m \le 1,\ \frac{1}{|\cell|}\int_\cell m   \,dx\le m_0\right\},
\end{equation}
for suitable constants $-\beta<m_0<0$ (in accordance with \eqref{eq:rem_case}).

Starting from early papers by Cantrell and Cosner \cite{MR1014659,MR1112065}, this question 
and related ones have been analyzed thoroughly in the literature, mainly in bounded domains, 
with different boundary conditions \cite{ly,llnp,MAZARI2022401}, or with different diffusion operators 
\cite{beroro,MR3771424,dipierro2021nonlocal,PePiSc}, also in the context of composite membranes 
\cite{MR1796024,MR2421158}. It turns out that the infimum of $\lambda_1(m)$ on $\Mcal$ is achieved 
by a bang-bang weight which saturates all the constraints. More precisely we have the following 
result.
\begin{theorem}\label{thm:existsbangbang}
Let $\lambda_1(m)$ and $\Mcal$ be defined as in \eqref{eq:weighted_eigenvalue} and \eqref{eq:Mcal}, 
respectively, and assume $-\beta<m_0<0$. Then
\[
\min_{m\in\Mcal} \lambda_1(m)
\]
is achieved, and for any minimizer $m^*$ there exists a open set 
$D^*\subset \cell$, with Lebesgue measure
\[
|D^*|=\delta := \frac{m_0+\beta}{1+\beta}|\cell|,
\]
such that
\begin{enumerate}
\item $m^*$ is a bang-bang weight:
\[
\left. m^*\right|_{\cell} = \ind{D^*} - \beta \ind{\cell\setminus D^*}
=
\begin{cases}
1  & \text{in }D^*\\
-\beta  & \text{in }\cell\setminus D^*,
\end{cases}
\]
and in particular $m^*$ has average equal to $m_0<0$;
\item $D^*$ is the superlevel set of the associated positive principal eigenfunction $u^*$:
\[
D^* = \left\{x\in \cell: u^*(x)> t^*\right\},
\]
for a suitable $t^*>0$.
\end{enumerate}
\end{theorem}
\begin{remark}\label{rmk:zeromeasurelevelsets}
In principle, since $m^*$ is $L^\infty$, $D^*$ is only a measurable set. 
On the other hand, using elliptic regularity and the 
equation, one can see that any optimal eigenfunction is of class $C^{1,\alpha}$, for every  
$0<\alpha<1$, and has negligible level sets, thus we can always choose $D^*$ to be an open set.
Actually, such regularity can be improved a lot, outside possible critical points of the eigenfunction, 
see Remark \ref{rmk:C11} ahead. 
\end{remark}

Statements analogous to Theorem \ref{thm:existsbangbang} have been proved in \cite{MR1014659}, partially, 
in the case of Dirichlet boundary conditions, and in \cite{ly} for Neumann ones, while the periodic case, 
for a related problem, 
was treated by Roques and Hamel in \cite{haro}. For these reasons we do not reproduce the 
proof here, but we rather refer to \cite[Appendix A]{haro} and \cite[Theorem 1.1]{ly}.

As a consequence of Theorem \ref{thm:existsbangbang}, the problem of determining the 
optimal spatial arrangement of resources for a species to survive is reduced to the following  
shape optimization problem with volume constraint, the unknown optimal shape $D^*$ 
being the favorable region of the habitat.

\begin{definition}\label{def:od}
For $\beta>0$ and $0<\delta <\frac{\beta}{1+\beta}|\cell|$, the optimal design problem for the survival threshold consists in solving
\[
\od(\delta) := \min\left\{\lambda(D):D\subset \cell,\ |D|=\delta \right\},
\]
where
\begin{equation*}
\lambda(D):=\lambda_1\left(\ind{D} - \beta \ind{\cell\setminus D}\right)
. 
\end{equation*}
\end{definition}
\begin{corollary}\label{coro:optdesprob}
For every $0<\delta <\frac{\beta}{1+\beta}|\cell|$ there exists (at least) a set $D_\delta\subset \cell$
and an associated eigenfunction $u_\delta\in W^{2,p}_{\per}(\cell)$, for every $1\le p<+\infty$, 
such that
\begin{enumerate}
\item $\od(\delta) = \lambda(D_\delta)$;
\item $u_\delta$ satisfies
\[
\begin{cases}
-\Delta u_\delta = \od(\delta) m_\delta u_\delta \quad\text{in }\cell\medskip\\
u_\delta\text{ is periodic, }u_\delta>0,\ \int_\cell u_\delta^2\,dx=1,
\end{cases}
\qquad\text{where }\left.m_\delta\right|_{\cell}:=\ind{D_\delta} - \beta \ind{\cell\setminus D_\delta};
\]
\item there exists $t_\delta>0$ such that
\[
D_\delta =  \left\{x\in \cell: u_\delta(x)> t_\delta\right\}.
\]
\end{enumerate}
\end{corollary}
After this discussion, the main questions concern the qualitative properties of $D_\delta$. Up 
to our knowledge, only few properties are known. 
It is a consequence of \cite{MR2214420,haro} (see also \cite{llnp}) that, up to a translation, 
$D_\delta$ is Steiner symmetric along each coordinate axis, that is, it is connected, symmetric with respect to all the 
hyperplanes $\{x_i = 0\}$ and convex along coordinate directions. On the other hand it is known that $D_\delta$ cannot be a ball \cite{llnp}.
Here we aim at providing a more detailed description of $D_\delta$, from both the qualitative and the 
quantitative point of view, in the singularly perturbed regime $\delta\to0$.

\subsection{Analysis of the small volume regime} \label{sec:mainres}

In a series of recent papers, in collaboration with Lorenzo Ferreri, Dario Mazzoleni and 
Benedetta Pellacci, we have discussed the above optimization problem, in the small volume regime 
$\delta\to0$, in a bounded domain $\Omega\subset\R^N$, 
with either Dirichlet or Neumann boundary conditions 
\cite{mapeve,ferreri_verzini,MPV3,ferreri_verzini2,FMPV24}. In these cases, the shape of the box 
$\Omega$ plays a very relevant role, and for general domains symmetrization techniques cannot be 
applied. On the other hand, in the small volume regime we expect concentration 
(and in particular connectedness) of $D_\delta$ at some point in $\overline{\Omega}$ and the main 
questions concern the asymptotic location of $D_\delta$, and its asymptotic shape. As one may expect, 
the two questions are intertwined. Roughly and informally, our results can be stated as 
follows.
\begin{itemize}
\item In the case of Dirichlet boundary conditions \cite{ferreri_verzini2}, if $\delta$ is 
sufficiently small, 
the optimal set $D_\delta$ is compactly contained in $\Omega$, connected, convex and nearly 
spherical (the definition of nearly spherical set is reported in Definition \ref{def:nearly_spherical} ahead). As $\delta\to0$, it concentrates at a point maximizing its distance to the lethal 
boundary, 
and the spherical asymmetry decays to zero in the $C^{1,1}$-norm; moreover, we can prove quantitative 
estimates of such decay, which is exponential, with respect to a negative power of $\delta$, 
in every $C^{1,\alpha}$-norm, $\alpha<1$. 
\item For Neumann boundary conditions \cite{MPV3,FMPV24}, $\overline{D}_\delta$ intersects 
$\partial\Omega$ for $\delta$ sufficiently small: more precisely, it is the (connected) 
intersection of $\Omega$ with a nearly spherical set, centered at a point of $\partial\Omega$. 
Now, as $\delta\to0$, concentration happens at a point of $\partial\Omega$ maximizing its mean 
curvature, while the spherical asymmetry decays in $C^{1,1}$ and only polynomially in 
$C^{1,\alpha}$. 
\end{itemize}

Here we want apply the same strategy to the periodic case. 
In this case, if one can prove concentration, then by translation invariance we can always assume that 
it happens at $0$. On the other hand, the study of the rate of decay of the spherical asymmetry 
requires a nontrivial combination of the techniques of the two cases above. Our main result is the 
following.
\begin{theorem}\label{thm:main}
Under the notations of Definition \ref{def:od} and Corollary \ref{coro:optdesprob}, 
there exists $\delta_{0}>0$ such that, for every $\delta\in (0,\delta_{0})$ and up to 
translations:
\begin{enumerate}
\item\label{pt:1} the restriction of $u_\delta$ to the period cell $\cell$ has a unique local maximum point, 
at $0$, and it is even in every variable $x_i$ and monotone decreasing  in $|x_i|$, for every $i=1,\dots,N$; 
\item\label{pt:2} $D_\delta$ is radially diffemorphic to a ball centered at $0$, i.e. there exists a function 
$\varphi_\delta\in C^{1,1}(\sphere^{N-1})$ such that 
\[
D_\delta=\left\{x: |x|< \delta^{1/N}\left(r_0+\varphi_\delta\left(
\frac{x}{|x|}\right)\right)\right\},
\]
where $r_0$ is such that $|B_{r_0}|=1$. Furthermore, $\varphi_\delta$ is even in every $x_i/|x|$.
\setcounter{fakecounter}{\theenumi}
\end{enumerate}
Moreover, as $\delta \to0$,
\begin{enumerate}
\addtocounter{enumi}{\value{fakecounter}}
\item\label{pt:3} there exists a positive constants $\Lambda_0$, $M$, depending on $N,\beta,\cell$, such that
\[
\od(\delta)=\lambda(D_\delta)=\delta^{-2/N}\left(\Lambda_0 +  O\left(\exp\left(-M\delta^{-1/N}\right)\right)\right);
\]
\item\label{pt:4} $\|\varphi_\delta\|_{C^{1,1}(\sphere^{N-1})}\to 0$;
\item\label{pt:5} for every $0<\alpha<1$ there exists a positive constant $M_\alpha$ such 
that 
\[
\| \varphi_{\delta} \|_{C^{1,\alpha}(\sphere^{N-1})} = O\left(\exp\left(-M_\alpha\delta^{-1/N}\right)\right).
\]
\end{enumerate}
\end{theorem}

Some comments are in order.
\begin{remark}\label{rmk:symmetries}
As we already mentioned, it has been shown in \cite{MR2214420,haro} (see also \cite{llnp}) 
that one can use consecutive Steiner symmetrizations along each coordinate axis. As a consequence 
we have that $D_\delta$ is Steiner symmetric (up to a translation, with respect to the origin): it is 
connected, symmetric with respect to all the hyperplanes $\{x_i = 0\}$ and convex along coordinate 
directions. Analogously,  $u_\delta$  has a maximum point at $0$, and it is even in every variable 
$x_i$ and monotone decreasing (in $\cell$) in $|x_i|$, for every $i=1,\dots,N$. In particular, both 
point \ref{pt:1} of Theorem \ref{thm:main}, and the symmetry properties of $\varphi_\delta$ in point \ref{pt:2}, 
are already well established.
\end{remark}
\begin{remark}\label{rmk:odzero}
Our strategy is inspired by the one used to study the mountain pass solutions of 
semilinear singularly perturbed problems, which was developed in many seminal papers 
at the end of last century, such as \cite{nitakagi_cpam,nitakagi_duke,NiWei:Spike95,MR1736974}.
In particular, the positive constant $\Lambda_0$ appearing in point \ref{pt:3} is related to a 
limit optimization problem on the entire space $\R^N$ which arises after a suitable blow-up 
procedure. We refer to Section \ref{sec:exp} for a description of the limit problem, which was 
already analyzed in \cite{MPV3}, and to Section \ref{sec:bu} for the blow-up argument. It is worth 
mentioning that, with respect to the semilinear case, several new difficulties arise, since  
the driving parameter $\delta = |D_\delta|$  is not explicit inside the equation and
the underlying equation is linear, with a discontinuous and
non-homogeneous weight, and with degenerate solutions (eigenfunctions).
\end{remark}
\begin{remark}\label{rmk:C11}
Being the optimal weight just $L^\infty$ and discontinuous, the blow-up procedure naturally provides 
just the $C^{1,\alpha}$ regularity of the eigenfunctions and of their regular level sets, for every  
$\alpha<1$. To recover the $C^{1,1}$ regularity we exploit the fact that the partial derivatives of 
the optimal eigenfunction satisfy a transmission problem, and we can use recent regularity 
results for such problems, obtained in \cite{Caffarelli2021:TransmissionProblems,Dong2021:TransmissionProblems}. Notice that such further 
regularity is not only interesting in itself, but it also allows to simplify drastically the 
quantitative part of our results, in particular point \ref{pt:5}.

Moreover, although the $C^{1,1}$ regularity is optimal at the level of the eigenfunctions, 
one expects the free boundary to be much more regular, since it is a regular level set of the 
solution of an elliptic equation. For instance, in the context of composite membranes, it is 
shown in \cite{MR2421158} that it is analytic, at least in dimension $N=2$.
\end{remark}
\begin{remark}\label{rmk:convex}
One significant implication of the $C^{1,1}$ decay in point \ref{pt:4} is that, as a consequence, 
$D_\delta$ is (connected and) strictly convex. 
\end{remark}
\begin{remark}\label{rmk:lamb}
It is a consequence of the results in \cite{llnp} that $D_\delta$ can not be a ball; in particular, we infer that 
$\varphi_\delta$ is nontrivial, for every $\delta$.
\end{remark}
\begin{remark}\label{rmk:neu_orthotope}
By Remark \ref{rmk:symmetries}, we have that the restrictions of the optimal weight $m_\delta$ and eigenfunction $u_\delta$ to 
\[
\Omega := \left(0,\frac{L_1}{2}\right) \times \dots \times \left(0,
\frac{L_N}{2}\right)
\] 
provide a minimizer of the optimal design problem on $\Omega$, with Neumann 
boundary conditions and prescribed measure $2^{-N}\delta$ (instead of $\delta$) (see also 
\cite[Section 2.7]{llnp}). Using Theorem \ref{thm:main} we can complement the results in 
\cite[Remark 1.7]{MPV3} obtaining that, also for the problem with Neumann boundary conditions 
in an orthotope, in the  small volume regime, the optimal set is nearly spherical, centered at 
a corner, it is convex and it meets the boundary orthogonally. Exponential decay of both the 
eigenvalues and of the spherical asymmetry holds true, too.
\end{remark}

The paper is structured as follows: in Section \ref{sec:exp} we recall the main properties of the 
limit problem on the entire space $\R^N$, and we exploit them to obtain a first asymptotic expansion 
of the optimal eigenvalue $\od(\delta)$; Section \ref{sec:bu} contains the core of our blow-up 
procedure, describing the convergence of a suitable rescaling of the optimal eigenfunctions and sets 
to those appearing in the limit problem; finally, in Section \ref{sec:quant_est} we present the 
quantitative estimates about the decay of the spherical asymmetry $\varphi_\delta$, together with the end of the proof of 
our main results.

\section{First order expansion of the optimal eigenvalue}\label{sec:exp}

In this section, to start with, we recall some properties of the rescaled limit eigenvalue problem 
in $\R^N$, which are crucial in our description. Next, we exploit such properties to provide a first 
estimate of the optimal eigenvalue $\od(\delta)$ as $\delta\to0$. 

For $A\subset\R^N$ (measurable and) bounded let us define the principal eigenvalue
\begin{equation}\label{eq:laRN}
\lambda(A, \R^N):= \inf\left\{\int_{\R^N}|\nabla w|^2 : 
w \in H^1(\R^N),\ 
\int_{A}w^2 - \beta\int_{\R^N\setminus A} w^2 = 1\right\}.
\end{equation}
Then, by symmetrization, one can see that the minimization of such eigenvalue, 
among sets $A$ having prescribed measure, is achieved by a ball, with an associated radial and 
radially decreasing eigenfunction. Precisely, we have that
\begin{equation}\label{eq:lambdazero}
\Lambda_0 := \inf_{|A|=1} \lambda(A, \R^N) = \lambda(B, \R^N),
\end{equation}
where $B\subset\R^N$ denotes the ball of measure one, centered at the origin, 
and $r_0$ denotes its radius:
\begin{equation}\label{eq:unitaryball}
B=B_{r_0},\qquad |B_{r_0}| =1.
\end{equation}
Actually, as in the periodic case, this minimization can be equivalently performed enlarging 
the class of weights from bang-bang ones to the more general family
\begin{equation}\label{eq:Mcal'}
\Mcal':=\left\{ m \in L^\infty (\R^N): 
-\beta \le m \le 1,\ \int_{\R^N} (m+\beta)\,dx\le (1 + \beta)\delta\right\},
\end{equation}
which is the natural version of \eqref{eq:Mcal} in this case.

The minimizer 
\begin{equation}\label{eq:lim_weight}
\widetilde m_0 := \ind{B}-\beta\ind{\R^N\setminus B}
\end{equation}
of both problems is unique up to translations, and in turn the optimal eigenvalue 
$\lambda(B, \R^N)$ is achieved by an eigenfunction $w\in H^1(\R^N)$, solution of 
\begin{equation}\label{eq:lim_prob}
-\Delta w = \Lambda_0 \widetilde{m}_0 w \quad \text{in } \R^N.
\end{equation}
which is positive, radially symmetric, radially decreasing and normalized in $L^2(\R^N)$ (and 
uniquely determined by such normalization). Actually, $w$ is explicit in terms of Bessel functions, and it decays 
exponentially at infinity:
\begin{equation}\label{eq:decay_w}
|w(x)|+|\nabla w(x)| \sim C |x|^{-(N-1)/2} e^{-\sqrt{\IM\beta}|x|}\qquad
\text{as }|x|\to+\infty;
\end{equation}
we refer to \cite[Section 2]{MPV3} for the proofs of these results, as well as for further 
details about $\Lambda_0$ and $w$.

Based on the previous discussion, we can state and prove the main result of this section.
\begin{proposition}\label{prop:firstordexp}
There exists a positive constant $C$ such that, as $\delta\to0$,  
\begin{enumerate}
\item $\od(\delta)\le \delta^{-2/N}\left(\Lambda_0 - C\exp(-2\sqrt{\IM\beta}\,\delta^{-1/N}) \right)$
\item $\od(\delta)\ge \delta^{-2/N}\left(\Lambda_0 + o(1)\right)$.
\end{enumerate}
\end{proposition}
\begin{proof}
The first part follows by considering the weight
\[
\hat m_\delta := \ind{\delta^{1/N}B}-\beta\ind{\cell\setminus \delta^{1/N}B},
\]
extended periodically, as a competitor in the definition of $\od(\delta)$, and consequently  
the restriction to $\cell$ of the function $w_\delta$, defined as
\[
w_\delta(x):= \delta^{-1/2}w(\delta^{-1/N} x),
\] 
as competitor in the definition of $\lambda(\delta^{1/N}B)$. We obtain, for $\delta$ sufficiently 
small (so that $\delta^{1/N}B\subset \cell$) and exploiting \eqref{eq:decay_w},
\[
\begin{split}
\od(\delta) &\le \frac{\int_{\cell}|\nabla w_\delta|^2}{\int_{\delta^{1/N}B}w_\delta^2 - \beta\int_{\cell\setminus\delta^{1/N}B} w_\delta^2} = 
\frac{\int_{\R^N}|\nabla w_\delta|^2 - \int_{\R^N\setminus\cell}|\nabla w_\delta|^2}{
\int_{\R^N}\widetilde m_0  w_\delta^2 + \beta\int_{\R^N\setminus\cell} w_\delta^2}\\
&\le \delta^{-2/N}\left(\Lambda_0 - C\int_{\R^N\setminus\delta^{-1/N}\cell}|x|^{-(N-1)} e^{-2\sqrt{\IM\beta}|x|}\right),
\end{split}
\]
whence the first part follows by direct calculations.

Concerning the second estimate, as we already noticed in Remark \ref{rmk:neu_orthotope}, we have that any optimal favorable set $D_\delta$ associated to 
$\od(\delta)$ is symmetric with respect to all the hyperplanes $\{x_i = 0\}$ (up to a translation). 
Then the restrictions of the corresponding weight $m_\delta$ and eigenfunction $u_\delta$ to 
\[
\Omega := \left(0,\frac{L_1}{2}\right) \times \dots \times \left(0,
\frac{L_N}{2}\right)
\] 
provide an admissible competitor (actually, they are a minimizer) for the optimal design problem on $\Omega$, with Neumann 
boundary conditions and prescribed measure $2^{-N}\delta$ (instead of $\delta$). Since 
\[
\frac{\int_{\Omega}|\nabla u_\delta|^2}{\int_{\Omega} m_\delta  u_\delta^2} =
\frac{\int_{\cell}|\nabla u_\delta|^2}{\int_{\cell} m_\delta  u_\delta^2} = \od(\delta),
\]
we can exploit the results in \cite{MPV3} (and in particular Remark 1.7 and equation (1.7) therein) to infer
\[
\od(\delta) \ge \frac14 \Lambda_0\cdot(2^{-N}\delta)^{-2/N} + o(\delta^{-2/N}),
\]
and also the second part follows.
\end{proof}

\section{Blow-up analysis}\label{sec:bu}

In order to refine the estimate from below in Proposition \ref{prop:firstordexp}, we resort to 
a blow-up analysis. (For a similar strategy in the Neumann or Dirichlet case, see 
\cite[Section 4]{MPV3} or \cite[Section2]{ferreri_verzini2}, respectively.) This will shed more light on the relation between the optimal design problem in 
Definition \ref{def:od} and the limiting one in \eqref{eq:lambdazero}, as depicted in Propositions 
\ref{prop:bu_w_conv} (for the optimal eigenfunctions and weights) and \ref{pro:nearlysphericalalpha} 
(for the optimal sets) below. Moreover, such blow-up analysis will provide uniform exponential 
decay of the rescaled eigenfunctions at infinity (Lemma \ref{lem:exp_decay_sol} and Corollary 
\ref{coro:exp_decay_grad}), which in turn will trigger the use of quantitative estimates in 
Section \ref{sec:quant_est}.

Recall that, after Corollary \ref{coro:optdesprob}, we denote with $u_{\delta}$ the 
positive eigenfunction associated to $\od(\delta)$, normalized in $L^2_\per(\cell)$, having 
associated optimal weight $m_\delta\in L^\infty_\per(\cell)$, with corresponding 
optimal favorable set $D_\delta\subset\cell$, with $|D_\delta|=\delta$. 

Moreover, as we already noticed in Remark \ref{rmk:symmetries}, after a suitable translation 
we can assume without loss of generality that $D_\delta$ is Steiner symmetric and that 
the restriction of $u_\delta$ to the period cell $\cell$ has a unique local maximum point, 
at $0$, and it is even in every variable $x_i$ and monotone decreasing  in $|x_i|$, for every $i=1,\dots,N$. In particular, we have that $u_\delta$ is radially decreasing (although not radially symmetric) in $\cell$.

We define the blow-up functions (centered at $0$, because of the previous discussion):
\begin{equation}\label{eq:bufuns}
\widetilde{u}_{\delta}(x) := \delta^{1/2} \, u_{\delta}
(\delta^{1/N}x), \qquad
\widetilde{m}_{\delta}(x) := m_{\delta}(\delta^{1/N}x),
\end{equation}
and sets
\begin{equation}\label{eq:busets}
\widetilde{\cell}_{\delta} := \delta^{-1/N}\cell,\qquad
\widetilde D_{\delta} := \delta^{-1/N} D_\delta =
\left\{ x \in \widetilde{\cell}_{\delta} : \widetilde{u}_{\delta}(x)>\widetilde t_\delta \right\},
\end{equation}
where $\widetilde t_\delta:= \delta^{1/2}\left. u_{\delta}\right|_{\partial D_\delta}>0$. 
Notice that with this choice we have, for every $\delta>0$,
\[
\int_{\widetilde\cell_\delta} \widetilde{u}_\delta^2 =| \widetilde D_\delta| = 1.
\]
Moreover, the functions $\widetilde{u}_{\delta}$ solve
\begin{equation}\label{eq:GeneralDiffProblemBUIndef}
\begin{cases}
-\Delta \widetilde{u}_{\delta} = \widetilde{\lambda}_{\delta} \widetilde{m}_{\delta} \widetilde{u}_{\delta} & \text{in } \widetilde{\cell}_{\delta} , \smallskip\\
\widetilde{u}_{\delta} \in H_{\text{per}}^1(\widetilde\cell_{\delta}),
\end{cases}
\end{equation}
where 
\begin{equation}\label{eqn:DefBlowUpEigenvIndef}
    \widetilde{\lambda}_{\delta}:=
    \delta^{2/N}\od(\delta)
\end{equation}
naturally inherits the variational characterization
\[
\widetilde{\lambda}_{\delta} = \min \left\{
\int_{\widetilde\cell_\delta} |\nabla u|^2 \,dx : A\subset\widetilde\cell_\delta,\ |A| =1,\ 
u \in H^1_\per(\widetilde\cell_\delta),\ \int_{\widetilde\cell_\delta} (\ind{A}-\beta\ind{\widetilde\cell_\delta\setminus A}) u^2  \,dx= 1\right\}.
\]
To start with we investigate the convergence of the blow-up sequence.
\begin{proposition}\label{prop:bu_w_conv}
Let $\Lambda_0$, $\widetilde m_0$ and $w$ be defined as in \eqref{eq:lambdazero}, \eqref{eq:lim_weight} 
and \eqref{eq:lim_prob}, respectively. Under the previous notation we have, as $\delta\to0$:
\begin{enumerate}
\item\label{pt:1b} $\widetilde\lambda_\delta \to \Lambda_0$;
\item\label{pt:2b} $\widetilde{m}_\delta\to \widetilde m_0$ weakly-$\ast$ in $L^{\infty}(\R^N)$;
\item\label{pt:3b} up to subsequences, there exists a constant $0<\sigma\le 1$ (depending on the subsequence
$\delta_n\to0$) 
such that 
\[
\widetilde u_\delta\to \sigma w
\] 
(strongly) both in $W^{2,p}_{\loc}(\R^N)$ and in 
$C^{1, \alpha}_{\loc}(\R^N)$, for every $1\le p <+\infty$ and $0 < \alpha < 1$. 
\end{enumerate}
\end{proposition}
\begin{proof}
The proof of \ref{pt:1b} is a direct consequence of Proposition \ref{prop:firstordexp} 
and the definition of $\widetilde\lambda_\delta$ in \eqref{eqn:DefBlowUpEigenvIndef}. In 
particular, in the following we can assume that $\delta$ is sufficiently small so that
\[
\widetilde\lambda_\delta \le \Lambda_0+1.
\] 

Concerning \ref{pt:2b} and \ref{pt:3b}, to avoid heavy notations, we keep writing $\widetilde{m}_\delta$, 
$\widetilde{u}_\delta$ and so on (with $0<\delta<\delta_0$, for a suitable $\delta_0>0$ small), 
although in principle we should take a sequence $\delta_n\to0$ as $n\to+\infty$, 
extract subsequences $(\delta_{n_k})_k$, and possibly use the Uryson subsequence principle, in 
case the limit is unique. We first notice that, since
\[
\|\widetilde m_\delta\|_{L^{\infty}(\R^N)} =\max\{ 1,\beta\}, 
\]
then there exists $m_\infty \in L^{\infty}(\R^N)$ such that, up to subsequences, 
\[
\widetilde{m}_\delta\to m_\infty
\qquad \text{weakly-$\ast$ in $L^{\infty}(\R^N)$.}
\]
Moreover we observe that $ m_\infty \in \Mcal'$, defined in \eqref{eq:Mcal'}, 
because such set is closed with respect to the weak-$*$ convergence. Now, let us fix 
any compact set $K\subset\R^N$. Then $K\subset \widetilde\cell_\delta$ for $\delta$ sufficiently small, 
and we can use equation \eqref{eq:GeneralDiffProblemBUIndef} and part \ref{pt:1b} to estimate
\[
\|\widetilde u_\delta\|_{H^1(K)}^2 \le \left(1 + \widetilde\lambda_\delta \|\widetilde m_\delta\|_{L^{\infty}
(\R^N)}\right) \|\widetilde u_\delta\|_{L^2(K)}^2 \le 1 + (\Lambda_0+1)\max\{ 1,\beta\}.
\]
Since this bound is independent of $K$ and $\delta$, we can consider a sequence $K_n \supset B_n(0)$ 
of compact sets and exploit a diagonal argument to infer the existence of $u_\infty\in H^1(\R^N)$ such that, 
up to further subsequences,
\[
\widetilde{u}_\delta\to u_\infty\text{ weakly in $H^1_{\loc}(\R^N)$},
\qquad\text{ and } -\Delta u_\infty = \Lambda_0 m_\infty u_\infty\quad \text{ (weakly) in } H^1(\R^N).
\]
Finally, one can apply standard elliptic regularity estimates to 
\eqref{eq:GeneralDiffProblemBUIndef}, 
based again on the uniform $L^\infty$ bounds of $\widetilde\lambda_\delta \widetilde m_\delta$, to bootstrap 
such convergence to the strong $W^{2,p}_{\loc}\cap C^{1, \alpha}_{\loc}$ one. In particular, 
$u_\infty$ is nonnegative, even in every variable $x_i$ and monotone decreasing  in $|x_i|$, for 
every $i=1,\dots,N$.

Thus, to conclude the 
proof, we are left to show that $m_\infty = \widetilde m_0$ and $u_\infty = \sigma w$, for some $0<\sigma \le1$. 
This will follow from the variational characterization of $\Lambda_0$, once we know that $u_\infty$ 
is not trivial.

To exclude that $u_\infty\equiv0$, we recall that, by construction, $\widetilde u_\delta$ is not 
constant and it is normalized in $L^2(\widetilde\cell_\delta)$. We obtain
\[
0< \int_{\widetilde \cell_\delta} |\nabla \widetilde u_\delta|^2 =\int_{\widetilde \cell_\delta} \widetilde m_\delta 
\widetilde u_\delta^2 = \int_{\widetilde D_\delta} \widetilde u_\delta^2 - \beta \int_{\widetilde \cell_\delta \setminus \widetilde D_\delta} \widetilde u_\delta^2 = (1+\beta) \int_{\widetilde D_\delta} \widetilde u_\delta^2 - \beta.
\]
Since $\widetilde u_\delta$ has a global maximum at $0$, we have
\begin{equation}\label{eq:intpositivo}
0<\frac{\beta}{\beta+1} < \int_{\widetilde D_\delta}\widetilde u_\delta^2 \le |\widetilde D_\delta| \cdot \widetilde u_\delta^2 (0)\to u_\infty ^2(0)
\end{equation}
by pointwise convergence, thus $u_\infty$ is nontrivial.

Once we know that $u_\infty$ is not trivial, we can test its equation with itself, obtaining
\[
0<\int_{\R^N} |\nabla  u_\infty|^2 = \Lambda_0 \int_{\R^N}  m_\infty   u_\infty^2,
\qquad\text{ i.e. }
\Lambda_0 = \frac{\int_{\R^N}|\nabla  u_\infty |^2}{\int_{\R^N}  m_\infty   u_\infty ^2} .
\]
Since $u_\infty\in H^1(\R^N)$, $m_\infty\in\Mcal'$ and $\int_{\R^N}  m_\infty   u_\infty ^2>0$, the 
variational characterization of $\Lambda_0$ (see the discussion in Section \ref{sec:exp}) implies 
that $u_\infty$ is a multiple of a translation of $w$ (and $m_\infty$ the corresponding translation 
of $\widetilde m_0$). Recalling the symmetry properties of both $u_\infty$ and $w$, we obtain that such 
translation is trivial, and $u_\infty = \sigma  w$, for some $\sigma \neq 0$. Finally, recalling that 
$u_\infty$ is non-negative, nontrivial, and with $L^2$ norm not exceeding $1$ 
(by strong $L^2_\loc$ convergence), we obtain $0<\sigma \le 1$, concluding the proof.
\end{proof}
\begin{remark}\label{rmk:notdefinedinRN}
Notice that, since neither $\widetilde u_\delta$ nor its restriction to $\widetilde \cell_\delta$ 
are in $H^1(\R^N)$, it is not proper to speak about global (either week or strong) 
$H^1(\R^N)$ convergence to $\sigma  w$, although  this will be the case after some suitable truncation. 
Actually, the same argument shows that one has always $\sigma =1$, for the full sequence $\delta\to0$.
See Remark \ref{rmk:strongconv} for more details. 
\end{remark}

The previous proposition has straightforward implications about the convergence of the rescaled 
optimal sets $\widetilde D_\delta$ to the unit volume ball $B=B_{r_0}(0)$, and more precisely of the
free boundary 
\begin{equation}\label{eq:impl_free_bd}
\partial\widetilde D_\delta = \left\{x \in \widetilde{\cell}_{\delta} : \widetilde{u}_{\delta}(x)=\widetilde t_\delta \right\}
\end{equation}
to $\partial B$. We start noticing that 
such optimal sets are uniformly bounded, from both outside and inside.
\begin{lemma}\label{lem:beta_outside}
For every $\overline{r}>r_0$ there exists $\delta_0>0$ such that, for every $0<\delta<\delta_0$, 
\[
\widetilde D_\delta \subset B_{\overline{r}}=B_{\overline{r}}(0),
\]
where $\widetilde D_\delta$ (and $\widetilde t_\delta$) are defined in \eqref{eq:busets}. In particular, under the same conditions,
\[
\widetilde m_\delta \equiv-\beta\qquad\text{in }\widetilde\cell_\delta\setminus B_{\overline{r}}.
\]
\end{lemma}
\begin{proof}
Fix $\overline{r}>r_0$, and assume by contradiction the existence of a sequence $\delta_n\to0$ such 
that (up to a relabeled subsequence)
\[
\widetilde D_{\delta_n} \setminus \overline{B}_{\overline{r}}\neq\emptyset\qquad\text{and}\qquad 
\widetilde u_{\delta_n} \to \sigma  w\text{ uniformly in } \overline{B}_{2\overline{r}}.
\]
Moreover, since $\widetilde D_{\delta_n}\ni0$ is open and connected, and 
$|\widetilde D_{\delta_n}|=|B|=1$, we have that $B\setminus \widetilde D_{\delta_n}
\neq\emptyset$. We deduce that both $\partial \widetilde D_{\delta_n} \cap \partial 
B_{\overline{r}}$ and $\partial \widetilde D_{\delta_n} \cap \partial B$ are non-empty, 
namely there exist sequences $(x_n)_n$, $(y_n)_n$ such that, for every $n$,
\[
|x_n|=r_0,\qquad |y_n|=\overline{r},\qquad \widetilde u_{\delta_n}(x_n) = 
\widetilde u_{\delta_n}(y_n) = \widetilde t_{\delta_n}
\]
(recall that $B=B_{r_0}$ denotes the ball of unit volume, and $r_0$ its radius). But this contradicts the local uniform convergence of $\widetilde u_{\delta_n}$ to $\sigma w$, since 
\[
\left.w\right|_{\partial B}>\left.w\right|_{\partial B_{\overline{r}}}.
\qedhere
\]
\end{proof}
\begin{lemma}\label{lem:1_inside}
For every $\underline{r}<r_0$ there exists $\delta_0>0$ such that, for every $0<\delta<\delta_0$, 
\[
\widetilde D_\delta \supset B_{\underline{r}}.
\]
In particular, under the same conditions,
\[
\widetilde m_\delta \equiv 1 \qquad\text{in } B_{\underline{r}}.
\]
\end{lemma}
\begin{proof}
The proof is analogous to that of Lemma \ref{lem:beta_outside}. Fix $\underline{r}<r_0$, and assume 
by contradiction that 
\[
\overline{B}_{\underline{r}} \setminus \widetilde D_{\delta_n}\neq\emptyset\qquad\text{and}\qquad 
\widetilde u_{\delta_n} \to \sigma  w\text{ uniformly in } \overline{B}_{2r_0},
\]
for some $\delta_n\to0$. Since $|\widetilde D_{\delta_n}|=|B|=1$, we have that 
$\widetilde D_{\delta_n}\setminus {B}_{r_0}\neq\emptyset$, and we can find  
sequences $(x_n)_n$, $(y_n)_n$ such that, for every $n$,
\[
|x_n|=r_0,\qquad |y_n|=\underline{r},\qquad \widetilde u_{\delta_n}(x_n) = \widetilde u_{\delta_n}(y_n) =
\widetilde t_{\delta_n},
\]
again contradicting the uniform convergence of $\widetilde u_{\delta_n}$ to $\sigma w$.
\end{proof}
\begin{corollary}\label{coro:inner_outer_balls}
We have, as $\delta\to0$, 
\begin{itemize}
\item $\widetilde t_\delta = \left. w\right|_{|x|=r_0}+o(1)$,
\item $(1-o(1))B\subset \widetilde D_\delta\subset (1+o(1))B$,
\item $\partial\widetilde{D}_{\delta}\subset (1+o(1))B\setminus (1-o(1))B$.
\end{itemize}
\end{corollary}

To proceed, we recall the definition of nearly spherical set. 
\begin{definition}\label{def:nearly_spherical}
A (open) set $A\subset\R^N$ is \emph{nearly spherical, of class $C^{k,\alpha}$}, 
if there exists $\varphi\in C^{k,\alpha}(\sphere^{N-1})$, with $\|\varphi\|_{L^\infty}\le r_0/2$, 
such that 
\[
A =\left\{x\in\R^N:|x| < r_0 + \varphi(\theta),\text{ for }\theta:=\frac{x}{|x|}
\in\sphere^{N-1}
\right\}
\]
(recall that $r_0$ is the radius of the ball of unit measure in $\R^N$: $B=B_{r_0}(0)$), 
and we denote $\varphi_A$ such $\varphi$.
\end{definition}

This notion was exploited in the Dirichlet or Neumann cases for a refined description of 
the free boundary, see \cite[Section 5]{ferreri_verzini2} or \cite[Section 4]{FMPV24}, respectively.
Here we have the following result.
\begin{proposition}\label{pro:nearlysphericalalpha}
For $\delta$ sufficiently small, $\widetilde{D}_{\delta}$ is nearly spherical of class 
$C^{1,1}$, parametrized by $\varphi_{\delta}$. In addition
\[
\|\varphi_{\delta}\|_{C^{1,1}(\mathbb{S}^{N-1})}\to 0, \qquad \text{as $\delta\to 0$}.
\]
\end{proposition}
\begin{proof}
The proof can be obtained along the lines of \cite[Propositions 3.11, 5.10]{ferreri_verzini2}. 
Actually, the same statement in $C^{1,\alpha}$, $\alpha<1$, instead of $C^{1,1}$, is a direct 
application of the implicit function theorem.
Indeed, in view of Corollary \ref{coro:inner_outer_balls}, we have that 
\[
\partial\widetilde{D}_{\delta}\subset  B_{\overline{r}}\setminus B_{\underline{r}},
\]
for some $0<\underline{r}<r_{0}<\overline{r}$. Taking into account \eqref{eq:impl_free_bd}
and using polar coordinates $\rho= r_0+\varphi =|x|\in [\underline{r},\overline{r}]$, 
$\theta=\frac{x}{|x|} \in {\mathbb S}^{N-1}$, we have that $\varphi_{\delta}$ is 
implicitly defined by
\[
F(\varphi,\theta):=\widetilde u_{\delta}((r_0+\varphi)\theta)=\widetilde t_\delta.
\] 
As $\widetilde u_{\delta}$ converges to $w$ in $C^{1,\alpha}(\overline{B}_{\overline{r}})$, 
$\alpha<1$, 
\[
\max_{\overline{B_{\overline{r}}\setminus B_{\underline{r}}}}\partial_{r} \widetilde u_{\delta}<\frac12\max_{\overline{B_{\overline{r}}\setminus B_{\underline{r}}}}\partial_{r} w <0,
\quad \text{for every $\theta \in {\mathbb S}^{N-1}$. }
\]
Then we can apply the implicit function theorem, obtaining $F(\rho,\theta)=t_{\delta}$ if and 
only if $\rho=\rho(\theta)$. This argument can be implemented locally near any 
$\theta_0 \in {\mathcal \mathbb{S}^{N-1}}$, so that by compactness we obtain a globally defined 
$\varphi_{\delta}(\theta)=\rho(\theta)-r_{0}$, of class $C^{1,\alpha}$, $\alpha<1$.
Furthermore, 
\[
\nabla \varphi_{\delta}=-\frac{\left(r_{0}+\varphi_{\delta}(\theta)\right)}{\partial_{\rho}F(\rho,\theta)}\nabla_{T}\widetilde u_{\delta}
=-\frac{\left(r_{0}+\varphi_{\delta}(\theta)\right) }{\partial_{\rho}F(\rho,\theta)}\left(
\nabla\widetilde u_{\delta}-(\nabla \widetilde u_{\delta} \cdot \theta)\theta\right).
\]
where $\nabla_{T}\widetilde u_{\delta}$ denotes the tangential component of the gradient of 
$\widetilde u_{\delta}$.
Then the $C^{1,\alpha}$ convergence of $\widetilde u_{\delta}$ to the radial function $w$ 
immediately yields, for every $ \alpha<1$,
\[
\|\varphi_{\delta}\|_{C^{1,\alpha}(\mathbb{S}^{N-1})}\to 0 \qquad \text{as $\delta\to 0$}.
\]

Now, since we are reduced to work in $B_{\overline{r}}\setminus B_{\underline{r}}$, with 
$C^{1,\alpha}$ convergence to $w$, the boundary conditions play no further role. As a 
consequence, the improvement of the previous estimate to $\alpha=1$ can be obtained exactly as 
in the case of Dirichlet boundary conditions, which is detailed in \cite[Section 5]{ferreri_verzini2}. 
The strategy is lengthy but now well established: the regularity comes from the observation that 
the derivatives of $\widetilde u_\delta$ satisfy a transmission problem, which enjoys the regularity 
properties obtained in \cite{Caffarelli2021:TransmissionProblems,Dong2021:TransmissionProblems}; 
the decay is based on a uniform validity of these properties, after a $C^{1,\alpha}$ diffeomorphism 
which maps $\partial \widetilde D_\delta$ into $\partial B$. All the details are discussed 
in \cite[Section 5]{ferreri_verzini2}.
\end{proof}

The final step in this section  consists in proving some uniform exponential decay of $\left.\widetilde 
u_{\delta}\right|_{\widetilde\cell_\delta}$ and 
its gradient at infinity. (A similar property was crucial in the analysis of the Neumann case, see 
\cite[Theorem 2.2]{FMPV24}, and in turn such result was inspired by a similar approach used 
in~\cite[Lemma~2.3]{delpino_flores} to estimate the best constant for Sobolev embeddings in bounded 
domains.)
\begin{lemma}\label{lem:exp_decay_sol}
There exist universal constants $C_1,C_2>0$ and $\delta_0>0$ such that,  for all $0<\delta< \delta_0$,
\begin{equation}\label{eq:expdecayudelta}
\widetilde u_\delta(x)\leq C_1  e^{-C_2|x|},\qquad  \text{for all }x\in \widetilde\cell_\delta.
\end{equation}
\end{lemma}
\begin{proof}
As a starting remark, we notice that for every $\eps>0$ there exists $R=R(\eps)>0$ such that, 
for $\delta$ sufficiently small, we have $0<\widetilde u_\delta(x)\leq  \eps$ for $|x|>R$, $x\in \widetilde\cell_\delta$. Indeed, 
choose $R$ in such a way that $w(x)=\eps/2$ for $|x|=R$. Then, if $\delta$ is sufficiently small, we 
have that $B_{R+1}=B_{R+1}(0)\subset\widetilde \cell_\delta$ and, by uniform convergence, $\widetilde u_\delta(x)\leq  
\eps$ for $|x|=R$. Then, being $\widetilde u_\delta$ radially decreasing in $\widetilde\cell_\delta$, 
the estimates propagates to all $|x|\in \widetilde\cell_\delta \setminus B_R$.

Then we claim that: \emph{There exist $R_0>0$ and $\nu_0>0$ such that for all $R>R_0$ and 
$\delta$ sufficiently small we have 
\begin{equation}\label{eq:iterxdecay}
\sup_{\widetilde\cell_\delta\setminus B_R}\widetilde u_\delta\geq 2\sup_{\widetilde\cell_\delta\setminus B_{R+\nu_0}}\widetilde u_\delta.
\end{equation}
}

Indeed, assume by contradiction that there exist sequences $R_n\to+\infty$, $\nu_n\to+\infty$, 
$\delta_n\to 0$ and $x_n\in \widetilde\cell_{\delta_n}$ such that $|x_n|\geq R_n+\nu_n$, 
\[
M_n:=\sup_{\widetilde\cell_{\delta_n}\setminus B_{R_n}}\widetilde u_{\delta_n}(x)
\qquad\text{and}\qquad
\widetilde u_{\delta_n}(x_n)=:\mu_n>\frac{M_n}{2}.
\]
Since each $\widetilde u_{\delta_n}$ is radially decreasing in $\widetilde\cell_\delta$, 
it is not restrictive to assume that 
\[
B_{\nu_n}(x_n) \subset \widetilde\cell_{\delta_n};
\]
indeed, it is enough to replace $\nu_n$, $x_n$ with $\tau\nu_n$, $(R_n+\tau\nu_n)x_n/|x_n|$ 
respectively, for a universal $0<\tau<1/2$ depending on $(L_1,\dots,L_N)$ (actually 
$\tau=1/(1+2^{1/N})$ works for any choice of $\cell$). 
Moreover, the starting remark implies that $\mu_n,M_n\to 0$ as $n\to+\infty$, 
while  Lemma \ref{lem:beta_outside} implies
\[
\widetilde m_{\delta_n} \equiv-\beta\qquad\text{ in }\widetilde\cell_{\delta_n}\setminus B_{R_n}.
\]
Defining the auxiliary functions
\[
v_n(y)=\frac{1}{\mu_n}\widetilde u_{\delta_n}(y+x_n),
\]
then $v_n(0)=1$, $0< v_n< 2$ if  $|y|<\nu_{n}$ (thus $|y+x_n|> R_n$)
and $v_n$ solves the equation
\[
-\Delta v_n=-\widetilde\lambda_{\delta_n}\beta v_n\qquad \text{ in }B_{\nu_n}(x_n).
\]
Then we can argue as in the proof of Proposition \ref{prop:bu_w_conv} to pass to the limit as 
$n\to+\infty$; we obtain that $v_n\to v$ locally uniformly (and locally $W^{2,p}$) to a positive solution of
\[
-\Delta v=-\Lambda_0 \beta v\qquad \text{in }\R^N,\qquad v(0)=1.\\
\]
We obtain a contradiction, since this equation does not admit bounded nonnegative, nontrivial 
entire solutions. Then the claim is proved and \eqref{eq:iterxdecay} is proved for every $R> R_0$.

To conclude, we iterate the claim \eqref{eq:iterxdecay}: taking without loss of generality 
$\nu_0>R_0$ we infer that, for every $k\in \N$ and $\delta$ sufficiently small, 
\begin{equation}\label{eq:diadicdecay}
\sup_{\widetilde\cell_\delta\setminus B_{k\nu_0}}\widetilde u_\delta(x)\leq 2^{-k}\sup_{\widetilde\cell_\delta\setminus B_{\nu_0}}\widetilde u_\delta(x)\leq 2^{-k}\widetilde u_\delta(0)
\leq 2^{-k} (w(0)+1).
\end{equation}
Now, for every $x\in \widetilde\cell_\delta$ (and thus $\delta^{1/N}<|x|\cdot\max\left\{\frac{L_1}2,\dots,\frac{L_N}2\right\}$), we can find $k\in \N$ such that 
\[
k\nu_0\leq |x|\leq (k+1)\nu_0,
\]
and \eqref{eq:diadicdecay} yields
\[
u_\delta(x)\leq C\cdot 2^{-k}\leq 2C\cdot 2^{-|x|/\nu_0},
\]
for $C=w(0)+1>0$, whence the lemma follows.
\end{proof}
Although this is not strictly necessary for the following, we can use the previous lemma 
to obtain decay estimates also for $\nabla \widetilde u_\delta$.
\begin{corollary}\label{coro:exp_decay_grad}
There exist universal constants $C_1,C_2>0$ and $\delta_0>0$ such that,  for all $0<\delta< 
\delta_0$,
\begin{equation}\label{eq:expdecaygradudelta}
|\nabla \widetilde u_\delta(x)|\leq C_1  e^{-C_2|x|},\qquad  \text{for all }x\in \widetilde\cell_\delta.
\end{equation}
\end{corollary}
\begin{proof}
The corollary follows from the previous lemma, elliptic regularity estimates, and the embedding of $W^{2,p}\subset C^{1,\alpha}$, 
for $p$ sufficiently large.
More precisely, let $\delta_0$ be as in the previous lemma, $0<\delta<\delta_0$, and $x$ be such 
that $B_{4}(x) \subset \widetilde\cell_\delta$. Then, using \eqref{eq:expdecayudelta}, we have that, 
\[
\|\widetilde u_\delta\|_{C^{1,\alpha}(B_{3}(x)\setminus B_{2}(x))}\leq C\|\widetilde u_\delta\|_{W^{2,p}(B_{3}(x)\setminus B_{2}(x))}\leq C \|\widetilde u_\delta\|_{L^p(B_{4}(x)\setminus B_{1}(x))}\leq C  e^{-C_2|x|}
\]
for a universal constant $C>0$. 
\end{proof}

\section{Quantitative estimates}\label{sec:quant_est}

In this section we turn to the quantitative estimates of the spherical asymmetry $\varphi_\delta$ and 
to the intertwined issue of refining the estimate from below of the optimal eigenvalue $\od(\delta)$, 
obtained in Section \ref{sec:exp}. This will allow to conclude the proof of our main results. 

The solution of both these issues is connected with the possibility to construct suitable competitors 
$v_\delta \in H^1(\R^N)$ for the limit problem \eqref{eq:lambdazero} for $\Lambda_0$, associated with the rescaled optimal weights and sets $\widetilde m_\delta$ and $\widetilde D_\delta$, and having weighted Rayleigh 
quotient related to $\widetilde \lambda_\delta$ (see the beginning of Section \ref{sec:bu} for the 
definition of the blow-up sequences). In turn, after the exponential decay estimates obtained in 
Lemma \ref{lem:exp_decay_sol}, to build $v_\delta$ it is enough to 
properly cut off the rescaled eigenfunctions $\widetilde u_\delta$.

For easier notation, let us assume that 
\[
\min\left\{\frac{L_1}2,\dots,\frac{L_N}2\right\}>1
\]
(of course, this is always possible after an initial suitable scaling), so that, in particular,
\[
B_{\delta^{-1/N}+1}=B_{\delta^{-1/N}+1}(0) \subset\subset \widetilde \cell_\delta
\]
for $\delta$ sufficiently small. Then, we choose a smooth cut-off function $\eta_\delta$ 
satisfying
\[
\eta_\delta := 
\begin{cases}
1 & \text{in }B_{\delta^{-1/N}}\\
0 & \text{in }\R^N \setminus B_{\delta^{-1/N}+1},
\end{cases}
\qquad 0\le \eta_\delta \le 1,
\qquad |\nabla\eta_\delta| \le 2\text{ in }\R^N.
\]
Finally, we define
\[
v_\delta = \eta_\delta \cdot \widetilde u_\delta \in H^1_0(\widetilde \cell_\delta) \subset H^1(\R^N).
\]
\begin{lemma}\label{lem:control_wdelta}
Under the above notation, there exist positive constants $\delta_0$, $M$, $C$ such that 
we have, for every $0<\delta<\delta_0$,
\[
\frac{\int_{\R^N}|\nabla v_\delta|^2}{
\int_{\R^N}\widetilde m_\delta  v_\delta^2} \le \widetilde \lambda_\delta + C e^{-M\delta^{-1/N}} .
\]
\end{lemma}
\begin{proof}
Testing by $\eta_\delta^2  \widetilde u_\delta$ the equation satisfied by $\widetilde u_\delta$ 
we obtain
\[
\begin{split}
\widetilde \lambda_\delta \int_{\R^N}\widetilde m_\delta  \eta_\delta^2  \widetilde u_\delta^2 &=  
\int_{\R^N}\nabla  \widetilde u_\delta \cdot \nabla(\eta_\delta^2  \widetilde u_\delta) = 
\int_{\R^N} \eta_\delta^2|\nabla  \widetilde u_\delta |^2 
+2 \eta_\delta \widetilde u_\delta \nabla \eta_\delta \cdot \nabla  \widetilde u_\delta \\
&= 
\int_{\R^N}  | \nabla(\eta_\delta  \widetilde u_\delta) |^2 - 
\int_{\R^N} \widetilde u_\delta^2 |\nabla\eta_\delta|^2, 
\end{split}
\]
that is,
\[
\frac{\int_{\R^N}|\nabla v_\delta|^2}{
\int_{\R^N}\widetilde m_\delta  v_\delta^2} \le \widetilde \lambda_\delta + 
\frac{\int_{\R^N} \widetilde u_\delta^2 |\nabla\eta_\delta|^2}{
\int_{\R^N}\widetilde m_\delta  v_\delta^2},
\]
and to conclude the proof we are left to estimate the (rightmost) remainder term. 

As long as its denominator is concerned, recalling that $\widetilde m_\delta<0$ where $0<\eta_\delta<1$, at least for $\delta$ sufficiently small, we have
\[
\int_{\R^N}\widetilde m_\delta  v_\delta^2 \ge \int_{\widetilde\cell_\delta}
\widetilde m_\delta  \widetilde u_\delta^2 = \frac{1}{\widetilde \lambda_\delta} 
\int_{\widetilde\cell_\delta} |\nabla \widetilde u_\delta|^2 \ge 
\frac{1}{\widetilde \lambda_\delta}\int_{B} |\nabla \widetilde u_\delta|^2 =
\frac{1}{\Lambda_0}\int_{B} |\nabla w|^2 + o(1)\ge C>0,
\]
as $\delta\to0$, by $H^1_{\loc}$ convergence of $\widetilde u_\delta$ to $w$; on the
other hand, we can estimate the numerator of the reminder term using Lemma \ref{lem:exp_decay_sol} 
(and the uniform Lipschitz bounds on $\eta_\delta$). We obtain
\[
\int_{\R^N} \widetilde u_\delta^2 |\nabla\eta_\delta|^2 \le 4C_1^2 \int_{B_{\delta^{-1/N}+1}\setminus B_{\delta^{-1/N}}} e^{-2C_2|x|},
\]
and the lemma follows.
\end{proof}

\begin{corollary}\label{coro:control_wdelta}
For every $\delta$ sufficiently small
\[
\Lambda_0 < \lambda(\widetilde D_\delta,\R^N) \le \widetilde \lambda_\delta + C e^{-M\delta^{-1/N}}. 
\]
\end{corollary}
\begin{proof}
It is enough to notice that, by the definition of $\lambda(A,\R^N)$ in equation \eqref{eq:laRN},
\[
\lambda(\widetilde D_\delta,\R^N) < \frac{\int_{\R^N}|\nabla v_\delta|^2}{
\int_{\R^N}\widetilde m_\delta  v_\delta^2}.\qedhere
\]
\end{proof}
\begin{remark}\label{rmk:strongconv}
As we mentioned in Remark \ref{rmk:notdefinedinRN}, by exploiting Lemma 
\ref{lem:exp_decay_sol}, Corollary \ref{coro:exp_decay_grad}, and a concentration-compactness 
argument, one can see that 
\[
v_\delta \to w\qquad\text{strongly in }H^1(\R^N),
\]
and that in Proposition \ref{prop:bu_w_conv} one has always $\sigma =1$, for the full sequence $\delta\to0$.
\end{remark}

In order to conclude the proof of our main results we need one last ingredient: the following 
sharp quantitative estimate related to problem \eqref{eq:lambdazero}, which was first introduced in 
\cite{ferreri_verzini2}. Earlier versions of stability inequalities  of this kind have already 
been established in \cite{CFL,Mazariquant}.

\begin{theorem}[{\cite[Theorem 1.4]{ferreri_verzini2}}]\label{thm:teoFV}
There exist positive constants $C,\eps$ such that, for all $C^{1,1}$ nearly spherical sets 
$A\subset\R^N$, centered at the origin and parametrized by $\varphi_A$, satisfying
\begin{enumerate}
\item $\bari(A)=0$,
\item $|A|=1$,
\item $\|\varphi_A\|_{C^{1,1}(\mathbb{S}^{N-1})}\leq \eps$,
\end{enumerate}
we have 
\[
\lambda(A,\R^{N})-\lambda(B,\R^{N})\geq C \|\varphi_A\|^2_{L^2(\mathbb{S}^{N-1})},
\]
where $\lambda(B,\R^{N})=\Lambda_0$.
\end{theorem}

We are finally ready to complete the proof of our main results.
\begin{proof}[End of the proof of Theorem \ref{thm:main}]
We already mentioned that point \ref{pt:1} and the symmetry properties in point \ref{pt:2} are already 
well established, see Remark \ref{rmk:symmetries}. The remaining part of point \ref{pt:2}, as well as 
point \ref{pt:4}, are indeed the content of Proposition \ref{pro:nearlysphericalalpha}. Turning to 
point \ref{pt:3}, the estimates from above was already contained in Proposition 
\ref{prop:firstordexp}, 
while the one from below follows by Corollary \ref{coro:control_wdelta}, recalling that, by 
definition,
\[
\widetilde \lambda_\delta = \delta^{2/N}\od(\delta).
\] 
Finally, also point \ref{pt:5} follows by Corollary \ref{coro:control_wdelta}. Indeed, as 
$\delta\to0$, we already know by point \ref{pt:4} that
\[
\|\varphi_\delta\|_{C^{1,1}(\sphere^{N-1})} = o(1);
\]
this fact, and the symmetry properties of $\widetilde D_\delta$, allow to apply Theorem \ref{thm:teoFV} to the rescaled optimal 
set, yielding
\[
\|\varphi_\delta\|_{L^2(\sphere^{N-1})} \le \lambda(\widetilde D_\delta,\R^N) - \Lambda_0 \le  C e^{-M\delta^{-1/N}},
\]
by Corollary \ref{coro:control_wdelta}. Using the Gagliardo-Nirenberg inequality to interpolate 
these two estimates, one can easily conclude that 
\[
\| \varphi_\delta \|_{C^{1,\alpha}(\sphere^{N-1})}\le C \|\varphi_\delta\|_{L^2(\sphere^{N-1})}^{(1-\alpha)/(4+N)}.
\]
The proof is completed.
\end{proof}

%
%
%
%
\bigskip

\textbf{Data Availability.} Data sharing not applicable to this article as no datasets were generated or analyzed during the current study.

\bigskip

\textbf{Disclosure statement.} The authors report there are no competing interests to declare.

\bigskip

\textbf{Acknowledgments.} Work partially supported by: PRIN-20227HX33Z ``Pattern formation in nonlinear 
phenomena'' - funded  by the European Union-Next Generation EU, Miss. 4-Comp. 1-CUP D53D23005690006.
The author is a member of the INdAM-GNAMPA group (``Gruppo Nazionale per l'Analisi Matematica, la Probabilit\`a e le loro Applicazioni -- Istituto Nazionale di Alta Matematica'').


\bibliographystyle{abbrv}
\bibliography{frac_eig.bib}

\medskip
\small
\begin{flushright}
\noindent 
\verb"gianmaria.verzini@polimi.it"\\
Dipartimento di Matematica, Politecnico di Milano\\ 
piazza Leonardo da Vinci 32, 20133 Milano (Italy)
\end{flushright}

\end{document}